\documentclass[11pt]{article}
\newcommand{\twist}{
\setlength{\unitlength}{1mm}
\mbox{
\begin{picture}(8,4)
\put(0,-3){\fig{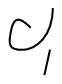}}
\end{picture}
}
}

\newcommand{\ltwist}{
\setlength{\unitlength}{1mm}
\mbox{
\begin{picture}(8,4)
\put(0,-3){\fig{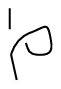}}
\end{picture}
}
}

\newcommand{\A}{\mbox {${\cal A}$}}
\newcommand{\1}{\mbox {${\bf 1}$}}

\newcommand{\fig}[1]{\includegraphics[scale=1]{#1}}

\newcommand{\newsec}[2]{\setcounter{section}{#1}
                         \setcounter{subsection}{0}
                         \large {\bf #1. #2 }\normalsize
                         \addcontentsline{toc}{section}{
                 \protect\numberline{#1}{#2}}
             }

\newtheorem{ppar}{}[section]
\newtheorem{lem}[ppar]{Lemma}
\newtheorem{theo}[ppar]{Theorem}
\newtheorem{prop}[ppar]{Proposition}
\newtheorem{corr}[ppar]{Corollary}

\newenvironment{ssec}{\begin{ppar} \rm}{\end{ppar}}

\usepackage{graphicx}

\begin{document}

\begin{center} \large
{\bf Relation between  quantum
invariants  of 3-manifolds and  2-dimensional
CW-complexes}\\
\normalsize
\vspace{0.7cm}
Ivelina Bobtcheva \\
{\em Department of mathematics, University of Ancona, Ancona, Italy}\\
Frank Quinn\\
{\em Department of mathematics, Virginia Polytechnic Institute and
State University, Blacksburg, VA, 24061, USA}
\thanks{Partially supported
through a grant from the US National Science Foundation.}\\
\vspace{0.7cm}
September, 2000
\end{center}

\begin{abstract}
We show that  the Reshetikhin-Turaev-Walker
invariant of 3-manifolds can be normalized to obtain an invariant of
   4-dimensional thickenings of 2-complexes.
  Moreover  when the underlying semisimple tortile category
  comes from the Lie family (``quantum groups'') over the ring
$Z_{(p)}[v]$ where $v$ is a primitive prime root of unity,
the 0-term in the Ohtsuki expansion of this invariant
  depends only on the spine and is
  the $Z/pZ$ invariant  of 2-complexes
  defined in \cite{Q:lectures}. As a consequence it is shown that
  when the Euler characteristic is greater or equal to 1,
  the 2-complex invariant depends only on homology.
  The last statement doesn't hold for the negative Euler characteristic case.
\end{abstract}

\newsec{1}{Introduction.}
The construction of ``quantum'' invariants of 3-manifolds was 
introduced by Reshetikhin-Turaev \cite{RT}
and put on a firmer footing by Walker \cite{W}, so we refer to these 
as ``RTW'' invariants.
For elaborate complete
  developments see the books
\cite{T} and \cite{KL}. In general the input for this construction is 
a finite semisimple tortile
category. Usually, however, the category is  understood to be one of 
the Lie family
obtained as  subquotients of representations of a deformation of the 
universal enveloping algebra of a
simple Lie algebra, specialize to a root of unity, see \cite{L:book}.

Analogous invariants of 2-dimensional CW complexes were introduced by 
the second author in
\cite{Q:lectures} and it was called in 
\cite{B2} the Q-invariant, so we will continue using this name. 
In general the input for this construction is a finite 
semisimple {\it symmetric\/} monoidal
category. Usually the category is taken to be one of the Lie family 
described by Gelfand-Kazhdan
\cite{GK}, obtained as subquotients of mod $p$ representations of 
simple Lie algebras. These invariants
are known to be invariant under deformations through 2-complexes. 
These deformations correspond to the
Andrews-Curtis moves \cite{AC} for presentations of groups. The 
(generalized) Andrews-Curtis conjecture
asserts that any simple homotopy equivalence of 2-complexes can be 
obtained by deformation through
2-complexes. This conjecture is expected to be false, and the Q- 
invariant was developed to try to
detect counterexamples.

The present work gives a connection between these invariants in the 
standard setting (input categories
  from the Lie family) by relating them both to an invariant of 
4-dimensional thickenings
of 2-complexes.  More explicitly we fix a simple Lie algebra and a 
prime $p$ greater than the dual
Coxeter number of the algebra. Let
$Z_{(p)}$ denote the integers localized at $p$ (= rationals with 
denominators prime to $p$), and $v$ be
a $p^{{th}}$ root of unity.

\begin{theo}
\label{thm} With Lie algebra and prime as above, there is $\hat 
Z(W)\in Z_{(p)}[v]$ defined for
$W$ a 4-dimensional thickening of a 2-complex such that
\begin{itemize}
\item[(a)]$\hat{Z}(W)$ is invariant under deformations through such
thickenings;
\item[(b)] there is a normalization of $\hat Z(W)$ in $Q[v]$  giving 
the RTW invariant of $\partial
W$; and
\item[(c)] the reduction mod $p$ of $\hat Z(W)$ is the Q- invariant 
of the spine of $W$.
\end{itemize}

  \end{theo}
  The normalization used in (b) and the proof of the theorem are given 
in \ref{mainproof}.

Part (c) provides topological interpretations for (at least some)   ``Ohtsuki
expansions'' of quantum invariants \cite{Oh}. In $Z_{(p)}[v]$ we can write
$$
\hat{Z}(W)=a_0+a_1(1-v)+\ldots +a_{p-1}(v-1)^{(p-1)}$$
with $a_i\in Z_{(p)}$. The mod $p$ reductions of the coefficients are 
well-defined, and are the Ohtsuki
expansion of $\hat Z(W)$.
$(1-v)$ is trivial in the mod $p$ reduction of $Z_{(p)}[v]$ so the 
mod $p$ reduction of $\hat Z$ is
equal to the mod $p$ reduction of $a_0$. In other words, the 0 term 
in the Ohtsuki
expansion of $\hat Z(W)$ is the Q- invariant of the spine of $W$.

The theorem also has implications for the Q- invariants.

\begin{corr}
\label{corra} When the Euler characteristic of the 2-complex 
is greater or equal to 1 its quantum invariant   is
determined by the homology of the complex. In particular, when the 
invariant 
is defined using a 
Gelfand-Kazhdan category it vanishes 
when the second
homology is nonzero, and otherwise is determined by the torsion 
subgroup of $H_1$.
\end{corr}
The corollary implies that the invariant cannot detect counterexamples to the 
original Andrews-Curtis
conjecture which concerns contractable complexes. This result was strongly
suggested by  numerical studies of the invariant (described at {\em
http://www.math.vt.edu/quantum\_topology}), but of course could not 
be proved that way.

In section 4 we use the corollary and computations for cyclic 
presentations to get an explicit formula
for invariants defined using the simplest Lie algebra. To state this 
recall $b_2$ (the second
Betti number) is the rank of $H_2(X;Z)$, and $t_1$ is the order of 
the torsion subgroup of $H_1(X;Z)$.

\begin{prop}\label{propa} Suppose $X$ is a 2-complex with Euler 
characteristic greater or equal to 1. Then
the class-0 $SL(2)$ Q- invariant of $X$ is $0$ if $b_2>0$ or if $p$ 
divides $t_1$, and is $t_1^{-2}\in
Z/pZ$ otherwise.
\end{prop}
Note the inverse in the second case is to be taken in $Z/pZ$. 
Explicitly, the invariant is the mod $p$
reduction of $r^2$, where $rt_1\equiv1$ mod $p$.

  The first author would like to thank Riccardo Piergallini for helpful
  discussions and Maria Grazia Messia for her encouragment and support.
\vspace{0.5cm}

\newsec{2}{Preliminaries and notations}

\begin{ssec}
Let $R$ be the
ring of integers localized at $p$, with a $p^{th}$ root of unity adjoined, so
$$R=Z_{(p)}[v]/<1+v+v^2+...+v^{(p-1)}>.$$
Let \A\ be a semisimple tortile category (see \cite{Shum}) over $R$
  and let
$S=\{1,a,b, \dots \}$ be a set of chosen representatives
for the simple objects in \A , i.e. $\1$ denotes the identity object
and small letters indicate simple objects. $A^{*}$ indicates the dual
object of $A$, and $A^{**}$ is canonically identified with $A$.
Furthermore, we use the
following notations:

\begin{eqnarray*}
&&\A \times \A \rightarrow \A \quad \mbox{(product)};\\
&\alpha _{A,B,C}:&(A B)
  C\rightarrow A (B C)\quad  \mbox{(associativity
  )};\\
&\gamma _{A,B} :&A B\rightarrow B A,
   \mbox{(commutativity)};\\
&\Lambda _{A}:&\1\rightarrow A^{*}\diamond A \mbox{ (coform)};\\
&\lambda _{A}:&A\diamond A^{*}\rightarrow \1 \mbox{ (form)};\\
&r _{A}=\lambda _{A}\Lambda _{A^*}: & \1\rightarrow \1 \mbox{ (rank)};\\
&\theta _{A}:&A\rightarrow A \,\1\stackrel{\Lambda _{A^*}}{\longrightarrow}
    A(A A^*)\rightarrow  (AA)A^*\\
    &&\stackrel{\gamma _{A,A}}{\longrightarrow}
    (AA)A^*\rightarrow A(AA^*)\stackrel{\lambda _{A}}{\longrightarrow}
    A\, \1
    \rightarrow A \mbox{ (twist)}\\
&dim(A,B)&\mbox{ is the dimension of $hom(A,B)$}.
\end{eqnarray*}

The main examples come from Lie algebras. Lustig \cite{L:book} 
defines a version of the quantum
enveloping algebra specialized at the root of unity $v$ as an algebra 
over $Z[v]$. Gelfand-Kazhdan
\cite{GK} define a subquotient category of the representations of 
this algebra. They show that when
reduced mod
$p$ the tensor product of representations induces a symmetric 
monoidal structure on the subquotient. It
follows from this that tensor product induces a monoidal structure 
(tortile, but no longer fully
symmetric) on the  subquotient with
$Z_{(p)}[v]$ coefficients. The point to check is that associativity 
maps are isomorphisms. But the hom
sets are finitely generated projective modules so mod $p$ isomorphism 
implies isomorphism over
$Z_{(p)}[v]$.   This conclusion can be considerably refined. General 
finiteness considerations show the
structure becomes monoidal over $Z[{1 \over Q},v]$ where $Q$ is some 
finite set of primes. We
conjecture it is sufficient to invert the primes less than the dual 
Coxeter number of the algebra.

The semisimplicity of the category implies
that given $a,b\in S$ there are bases for spaces of morphisms to or 
from the product. For every
simple object $z$ choose bases
$inj _k(z,ab):z\rightarrow ab$, and $proj _k(ab,z):ab\rightarrow z$,
$1\leq k\leq dim(z, ab)$, such that:

\begin{eqnarray*}
&& proj _l(ab,z)\circ inj_k(z,ab)=\delta _{k,l} id_z; \\
&& \sum _{z,k}inj_k(z,ab)\circ proj _k(ab,z) = id_{ab}.
\end{eqnarray*}

Using the coherence results in
\cite {Shum} we can represent certain morphisms in \A\ by labeled 
tangle diagrams (with
blackboard framing).
  Let $\underbar{A}=
(A_{1}, \,A_{2}\,\ldots ,A_{k})$
and  $\underbar{B}=(B_{1}, \,B_{2}\,\ldots ,B_{l})$ be sequences of objects
in \A .
Chose bracketings ${\cal B}_1$ and ${\cal
B}_2$ of the sequences to specify a way to form the products. Then
there is a one to one correspondence between the set of equivalence
classes of labeled
tangle diagrams $D:  \underbar{A}\rightarrow \underbar{B}$ and the
set of morphisms ${\cal B}_1(\underbar{A})\rightarrow {\cal
B}_2(\underbar{B})$ in \A\ formed by composition of products of
the elementary morphisms $\alpha , \, \gamma,\, \lambda ,\, \Lambda$
and  identities. The correspondence is
determined by the images of the elementary morphisms as
shown in figure \ref{elem.diagr}.

\begin{figure}[h]
\setlength{\unitlength}{1cm}
\begin{center}
\begin{picture}(10,7)
\put(0,0){\fig{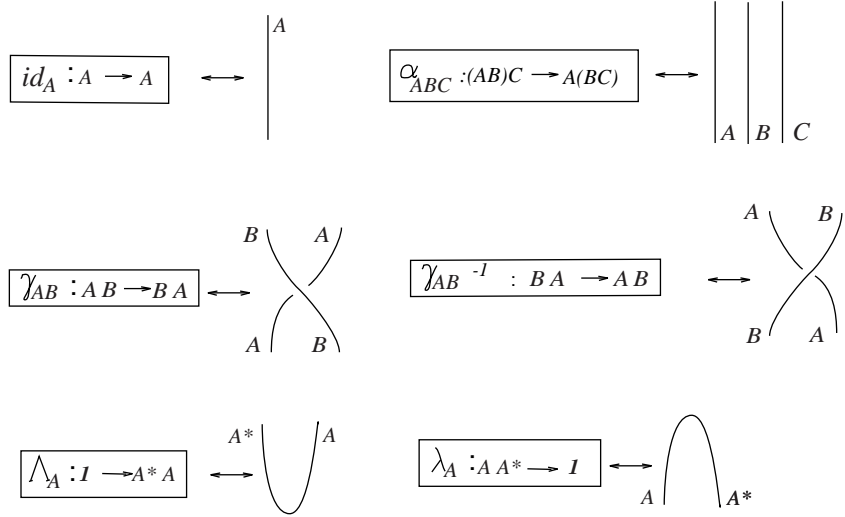}}
\end{picture}
\end{center}
\caption{Elementary diagrams.}
\label{elem.diagr}
\end{figure}

The equivalence relation on diagrams is generated by isotopy in
the plane and modified Reidemeiser moves, where the first  move
(involving the twist) is replaced by $$\twist =
\ltwist .$$

We will use labels that are ``linear combinations'' of objects with 
coefficients in $R$, and
associate morphisms to these as follows. The morphisms defined with 
genuine object labels are additive
in each label, so
$$\mbox{morph}(\cdots, A+B,\cdots)=\mbox{morph}(\cdots,
A,\cdots)+\mbox{morph}(\cdots, B,\cdots).$$
  Therefore if $r, s\in R$ we can define
$$\mbox{morph}(\cdots, rA+sB,\cdots)=r\,\mbox{morph}(\cdots, 
A,\cdots)+s\,\mbox{morph}(\cdots,
B,\cdots).$$
In particular we use the ``universal'' label $U=\sum _{a\in S}r_a a$ 
where simple objects are weighted
by their ranks. This is different from the normalization used in 
other constructions, and is done to
keep the constructions in the ring $R$ as long as possible.
\end{ssec}

\begin{ssec}
\label{def.nond}
A semisimple tortile category \A\ is called {\it nondegenerate} if the
following three conditions are satisfied:
\begin{itemize}
\item[(a)] $X^2=\sum _{a\in\Sigma}r_a{}^2$ is not a zero-divisor;
\item[(b)] $C_+=\sum _{a\in\Sigma}r_a{}^2\theta _a$ and
            $C_-=\sum _{a\in\Sigma}r_a{}^2\theta _a{}^{-1}$ are not 
zero-divisors;
\item[(c)] For any $a\in\Sigma$, $a\neq \1 $ , there exist
labeled tangle diagrams  $D$ and $D'$ such that each of them contains a
closed tangle component labeled by $a$.
Moreover, $D$ and $D'$ are equal if this
closed component is deleted, but the morphisms in \A\ corresponding to
$D$ and $D'$ are different.
\end{itemize}

The constructions of \cite{KL, RT, T, W} give 3-manifold invariants 
from nondegenerate tortile
categories. These invariants take values in the ring with inverses 
for the elements in (a) and (b)
adjoined, not necessarily in the original ring.

  The Lie family of categories are  nondegenerate in this sense.
These categories also satisfy:
\begin{itemize}
\item[(d)] $\theta _{a}$ for a simple
object $a$ acts as  $v^{t(a)} id_a$ for some
integer $t(a)$.
\end{itemize}
 From now on we
assume that \A\ is a semisimple tortile category over $R=Z_{(p)}[v]$
satisfying the conditions (a)--(d).
\end{ssec}

\begin{ssec}
\label{prop.nond}
We list two properties of nondegenerate
categories, which show how special such categories are (for the proofs
see \cite{W}).

The first property states the invariance under a
band-connected sum or difference of one tangle components with
another closed component. Let a tangle diagram
$D$ contain a closed component $K$ labeled with the
universal label $U=\sum _{a\in S}r_a a$, and $D'$ is obtained from $D$ by
  sliding any other tangle component along
$K$.
Then if \A\ satisfies (a) and (b), the morphisms
corresponding to the original
tangle diagram and the new one are the same.

The second property requires that the three conditions (a)--(c)
are
satisfied and describes the morphism corresponding to the diagram
$$T(b_1,b_2\dots ,b_k):(b_1,b_2\dots ,b_k)\rightarrow (b_1,b_2\dots
,b_k)$$ consisting from $k$ straight segments labeled by the $b_{i}$'s
and a unknoted closed component  which goes around them, labeled
by $U$. The  diagram is shown in figure \ref{proj.oper}.
It is shown in \cite{W} that
$$
T(b)=\left\{\begin{array}{cc} X^2 \,id_1 , &\mbox{if $b=\1$};\\
                                0 , &\mbox{ if $b\neq \1$}.
\end{array} \right.
$$
Then from the semisimplicity of the category it follows that
$$
T(b_1,b_2\dots ,b_k)=X^2\sum _{i=1}^{dim(\1 ,b_1b_2\dots b_k)}
inj_i(1,b_1b_2\dots b_k)\circ proj_i(b_1 b_2\dots b_k, \1 ).
$$
This shows $\frac{1}{X^2}T(b_1,b_2\dots ,b_k)$ is a morphism in the
category (i.e.\ over $R$, in spite of division by $X^2$).  More 
precisely it is the trace of the
projection over the trivial summand.

\begin{figure}[h]
\setlength{\unitlength}{1cm}
\begin{center}
\begin{picture}(4,2.5)
\put(0,0){\fig{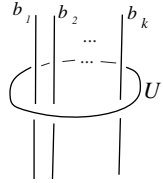}}
\end{picture}
\end{center}
\caption{The projection morphism $T(b_1,b_2\dots ,b_k)$.}
\label{proj.oper}
\end{figure}

 From the evaluation of $T$ above after
doing a handle slide, it follows that
$$
C_{+}C_-=X^{2}.
$$
\end{ssec}

\begin{ssec}
\label{def.phip}
Let
$\phi _p : R\rightarrow Z_p$ denote reduction mod $p$. Explicitly 
this is given by
$$
\phi_p( a_0+a_1 (1-v)+ \dots a_{l}(1-v)^{l})= a_0 \mbox{ mod $p$}.
$$
The following explains the role of hypothesis  \ref{def.nond}(d) in 
the relationship between the Lie
type tortile categories over $R$ and the Gelfand-Kazhdan categories 
over $Z/pZ$.

\begin{prop}
Let \A\ be a semisimple tortile category over $R$  satisfying 
condition \ref{def.nond}(d) and let
$\phi_p(\A )$ be the  category over $Z/pZ$ having the same objects as 
\A, and morphism modules the mod
$p$ reductions of morphisms in \A. Then $\phi_p(\A )$ is finite 
semisimple and symmetric.
\end{prop}

It is clear that mod $p$
reduction of structure in \A\ gives a finite semisimple tortile 
category. The twists and
commutativities satisfy
$$
\theta _{ab}=\theta _a^{-1}\theta _b^{-1} \gamma _{b,a}\gamma_{a,b}.
$$
For the quotient to be symmetric we need  $\gamma 
_{b,a}\equiv\gamma_{a,b}^{-1}$ mod $p$, or more
precisely
  the mod $p$ reductions of the twists $\theta_A$ should all be identities.
Since $\theta$ is additive it is sufficient to show this for simple objects. 
Hypothesis
\ref{def.nond}(d) asserts that these are powers of $v$ times identities.
$\phi_p$ takes powers of $v$ to 1, so the mod $p$ reductions are 
identity morphisms as required.
\end{ssec}

\newsec{3}{Proof of the Main Result.}
\label{goingup}
\begin{ssec} {\bf RTW Invariants}\quad
Let
$M$ be a closed  3-manifold, obtained via surgery on a framed
link $L$. Define $Z(L)$ to be
the morphism $\1\rightarrow \1$ in \A\  induced by the diagram
with underlying framed link $L$ and label $U$
on every component.
Then the quantum invariant of $M$
  is the element of $R[\frac{1}{X},\frac{1}{C_-},\frac{1}{C_+}]$ defined by
$$
Z_{RTW}(M)=\frac{1}{C_+^{\sigma_{+}}C_-^{\sigma_{-}}X^{\sigma_{0}}} Z(L)
$$
where $\sigma_{+}$, $\sigma_{-}$ and
$ \sigma_{0} $ are the
numbers of positive,
negative and  zero eigenvalues of the linking matrix of $L$.

The unnormalized  $Z(L)$ is a link invariant, and is
invariant under the second Kirby move (band-connected sum or
difference of two link components). It does not define a manifold 
invariant because it is not invariant
under the first Kirby move (adding or deleting an unknot of framing 
$\pm 1$ away of the
rest of the link). The normalization
factor in  $Z_{RTW}(M)$ compensates for this non-invariance, so 
$Z_{RTW}(M)$ is a manifold invariant.
Note the normalized invariant takes values in a larger ring because 
(in our cases) $C_{\pm}$ and $X$
are not invertible in $R$.
\end{ssec}

\begin{ssec} {\bf 4-thickening invariants}\quad
\label{def.4tick}
A  4-thickening $W$ of a 2 CW-complex $P$ is a 4-dimensional manifold 
together with a
  decomposition as a handlebody with 0-, 1-, and 2-handles
and an identification of the spine of the handlebody structure with 
$P$. Two 4-thickenings will
be called 2-equivalent if they can be deformed into each other
by a 2-deformation, i.e. by 1- and 2-handle slides and 0-1 and 1-2 handle
cancellations or introductions. We observe that a 2-equivalence of 
4-thickenings induces
a 2-equivalence (in the 2-complex sense) of their spines.

A 4-thickening with a single 0-handle can be described by a framed link
in $S^{3}$, obtained from the attaching maps of the 1-
and 2-handles.
The 1- handles are being  represented by
dotted unknots of framing 0 (the unknot represents the attaching
map of a canceling 2-handle), and the 2- handles correspond to the
undotted components.  An example (using the  blackboard framing) is 
shown in figure \ref{example}.
  Then (\cite{K:book}) two 4-thickenings are 2-equivalent if and
only if the
corresponding framed links can be deformed into each other through:
\begin{itemize}
\item[(a)] isotopy of framed links;
\item[(b)] handle moves
\subitem i) band-connected sum or difference of two dotted link components
(sliding an 1-handle over another 1-handle);
\subitem ii) band-connected sum or difference of two
undotted link components (sliding a 2-handle over a 2-handle);
\subitem iii) band-connected sum or difference of one undotted link 
component with
one dotted link component (sliding a 2-handle over 1-handle);

\item[(c)] any pair of one  dotted component $C$ and  one undotted
component $D$ can be removed or added if the geometric intersection
number of $D$ and the Seifert surface of $C$ is $\pm 1$ (1-2 handle
cancellation or introduction).
\end{itemize}

\begin{figure}[h]
\setlength{\unitlength}{1cm}
\begin{center}
\begin{picture}(8,5.5)
\put(0,0){\fig{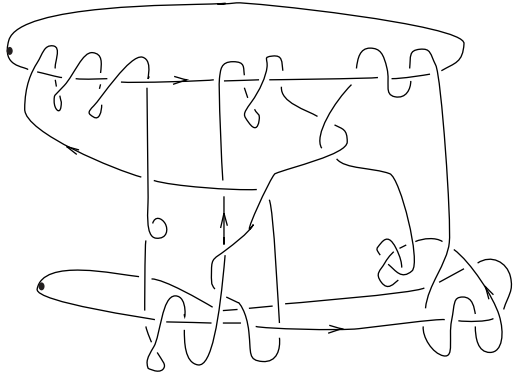}}
\end{picture}
\end{center}
\caption{A 4-thickening with spine described by the presentation
\newline
$P=<x,y|x^3y^2x^2y^{-1}, x^{-2}y^{2}>$.}
\label{example}
\end{figure}

Suppose  $L$ is a framed link diagram describing a 4-thickening. The 
link invariant
$Z(L)$ is unchanged by moves (a) and (b) but is not
invariant under (c). For one thing we have not distinguished between 
dotted and undotted
components, which amounts to exchanging all 1-handles with their 
cancelling 2-handles. The RTW
normalization accepts this exchange and normalizes $Z(L)$ to be invariant under
the first Kirby move. Here we use a different normalization that 
records something about the
1-handles.

\begin{prop} Let $W$ be a 4-thickening represented by a framed link
$L$ with $n$ dotted components. Then
$$
\hat{Z}(W)=\frac{1}{X^{2n}} Z(L)
$$
is  invariant under a 2-deformation.
\end{prop}

The only thing  to check is the invariance under (c).
Move (c) corresponds to changing $L$ into $L'$ by introducing
a dotted unknot $C$  and an undotted link  component
  $K$ such that the geometric intersection  number of $K$ and the disk
  bound by $C$ is $\pm 1$.
  According to the second property of
nondegenerate categories in \ref{prop.nond}, $C$ can be removed, and
the label of $K$
changed to the trivial label. But a component with the trivial label
doesn't contribute to $\hat Z(L')$, i.e.
$$
Z(L')= X^{2}\; Z(L).
$$
\end{ssec}

\begin{ssec} {\bf Proof of \ref{thm}}\quad
\label{mainproof}
Let \A\ be a  tortile category over $R$ satisfying the nondegeneracy 
conditions \ref{def.nond}(a)--(d),
and let $W$ be a 4-thickening of a 2-complex
$P$. Then to prove \ref{thm} we  show
\begin{itemize}
\item[(a)]$\hat{Z}(W) \in R$;
\item[(b)] the mod $p$ reduction
$
\phi _p(\hat{Z}(W))=Z_{Q}(P)$ is a 2-deformation invariant of the spine $P$.
\end{itemize}

\noindent{\bf Proof of (a).}
  Let $L$ be a link  describing $W$, and let $n$ be the number of 
dotted components of $L$.
$\hat{Z}(W)=\frac{1}{X^{2n}}Z(L)
$
and $Z(L)\in R$, so the only problem  comes from the factor
$\frac{1}{X^{2n}}$ (recall $X$ is not invertible in $R$).
But the diagram for $L$ can be deformed into a composition of
tangle diagrams, such that $n$ of the composition factors are of the form
$L'_k\otimes T(a_{i_{1}},a_{i_{2}},\dots ,a_{i_{s(k)}})\otimes L''_{k}$,
where $1\leq k\leq n$,
and $s(k)$ is the sum of the absolute values of all the exponents of
the $k$'th generator in all relations. Here we are regarding the 
2-complex $P$ as a presentation, with
generators the 1-cells and relations the 2-cells. Recall (\ref 
{prop.nond}) that
  $\frac{1}{X^{2}}T(a_{i_{1}},a_{i_{2}},\dots
,a_{i_{s(k)}})$ is in $R$. Since there are $n$ of these factors we 
see that we can divide $Z(L)$ by
$X^{2n}$ in $R$.

\medskip
\noindent{\bf Proof of (b).}
To prove (b) it is sufficient to identify the mod $p$ reduction as 
the Q- invariant, since this is
known to be a 2-complex invariant. However we use a different 
approach that also gives a proof of
\ref{corra}. In this approach we define $Z_{Q}$ by lifting: applying 
the 4-thickening invariant to a
canonical thickening of the spine, then reducing mod $p$.

Suppose
$P$ is a presentation. This means $P$ is a 2-complex with a single 
0-cell, an order and orientations for
the 1-cells, and attaching maps for the 2-cells expressed as words in 
the 1-cells. Then there is a
standard thickening $W_P$ defined by Huck
\cite{H}. This thickening is described by a framed link $L_{P}$ with 
$n+m$ components, such that any
component is an unknot of framing 0 and the geometric intersection
number of any two undotted or any two dotted components is 0.
The precise definition of $L_{0}$ can be given given as a closure of a braid.
Let $B_{n+m}$ be the braid group on $n+m$ strings. Let $y_{j}$,
$j=2\ldots n+m$ be the
braid group generator corresponding to interchanging the places of
the $j-1$ and the $j$ string. Then let
$r_{j,k}$, $j< k$ be the group element which moves the $j$ string to the $k$
place:
$$
r_{j,k}=y_{j+1}y_{j+2}\ldots y_{k-1}y_{k}.
$$
For each $j=1\ldots m$ we define a homomorphism $\psi _{j}$ from the free group
$F_n$ on
$n$ generators $x_1,x_2,\dots ,x_n$ into $B_{n+m}$ such that
$$
\psi _{j}(x_{k})=r_{j,k-1} y_{k+m}^2 (r_{j,k-1})^{-1},
$$
where we assume $r_{ii}=1$.
The link $L_P$ is defined to be the closure of the braid:
$$
\hat{P}=\psi _{1}(R_1) \psi _{2}(R_2) \ldots \psi _{m}(R_m).
$$
An example is given in figure \ref{braid}.

\begin{figure}[h]
\setlength{\unitlength}{1cm}
\begin{center}
\begin{picture}(3,7.5)
\put(0,0){\fig{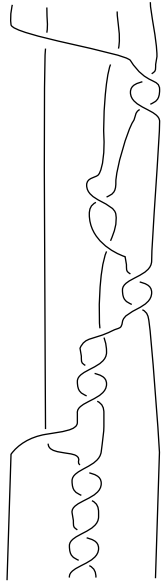}}
\end{picture}
\end{center}
\caption{The braid $\hat{P}$ for the 2-complex $P=<x,y|x^2, xyx^{-1}y^{-1}>$.}
\label{braid}
\end{figure}

For the purposes here we define $Z_{Q}(P)=\phi _{p}(\hat Z(W_P))$.

Next
we show this is the mod $p$ reduction of the general $\hat Z$ 
invariant. First observe that if $L$ is
the link describing an arbitrary thickening
$W$ of
$P$, its diagram  can be changed into the one for $W_P$ by (1) a 
finite number of
crossings  not involving dotted components; (2) adding twists
on undotted components; and (3) Reidemeister's moves. But it follows from
from \ref {def.phip} it follows that $ \phi _{p}(\hat Z(W))$ is
unchanged under these operations, i.e., $ \phi _{p}(\hat Z(W))= \phi 
_{p}(\hat Z(W_P))$.
This is true  because
the diagram can be sliced so that
  the changes are performed in  slices
not containing a dotted component.
Crossing involving a dotted component must be avoided because this 
takes the invariant out of the ring
$R$.

Finally we show that $Z_{Q}(P)$, as defined here, is invariant under 
a 2-deformation.   For this it is
enough to see that for a single 2-deformation move
$f:P\rightarrow P'$ there exist thickenings $W$ of $P$ and $W'$ of
$P'$, and a 2-deformation of thickenings $F: W\rightarrow W'$ that 
realizes the move on the spines.
The various sliding moves, and addition of a cancelling 1-2 pair of 
handles can be done starting with
any $W$, and we let $W'$ be the result. To remove a cancelling pair 
of handles we begin with some $W'$
and let $W$ be the thickening obtained by introducing handles.
\end{ssec}

\begin{ssec}
\label{r.g.change}{\bf Proof of \ref{corra}}\quad
Suppose $P$ is a presentation. Let $L$ be the framed link constructed 
above  to define the thickening
$W_P$. $L$ has the feature that the undotted components form a 
trivial framed link. There is
a ``duality'' defined for such framed links by interchanging the 
dotted and undotted sublinks.
Geometrically this corresponds to replacing 1-handles by their 
cancelling 2-handles, and replacing
2-handles with a particular set of cancelling 1-handles. Denote the 
spine of the associated handlebody
by
$P^*$, and the accisiated handlebody with
$(W_P)^*$. Huck
\cite{H} observes that there is an embedding of $W_P$ in $S^4$ with 
complement $(W_P)^*$.

If $P$ is given in terms of generators and relations by
$$P=<x_1,x_2,\dots ,x_n|\; R_1,R_2,\dots ,R_m>$$
then the
  dual presentation can be described
explicitly as
  $$P^*=<r_1,r_2,\dots ,r_m|\;
X_1,X_2,\dots ,X_n>.$$ The relations are
$$
X_k=r_{1}^{f_k^1} r_{2}^{f_k^2}\dots r_{m}^{f_k^m}.
$$
where $f_k^l$ is the total exponent of $x_k$ in $R_l$. Recall that 
the total exponents are entries in
the boundary homomorphism in the cellular chain complex. Since the 
chain complex is characterized up to
deformation by the homology, it follows that the dual presentation 
(up to a 2-deformation)
depends only on the homology of $P$.

Recall that $\hat Z$ is defined by normalizing a link invariant. Specifically
$$\hat Z(W) = \frac{1}{X^{2n}}Z(L)$$
where $n$ is the number of 1-handles (dotted components) of $W$ and 
$L$ is the link obtained by
forgetting dots. Since $W_P$ and $(W_P)^*$ have the same underlying 
link, and $n$ and $m$ 1-handles
respectively, we get
$$X^{2n}\,\hat Z(W_P) =X^{2m}\,\hat Z((W_P)^*).$$
In particular, if $m\geq n$, $Z_Q(P)=\phi _p(X^{2(m-n)})Z_Q(P^*)$.
Since $P^{*}$ depends only on the homology of $P$, it
follows that in this case, the invariant depends only on homology.


To complete the proof of the corollary we need to refine this in 
terms of Betti numbers and torsion in the case that the invariant is 
defined using a Gelfand-Kazhdan category. The important point is that 
for such categories $\phi _p(X^{2})=0$.

Recall  $b_i$ is defined by $H_i(P;Z) \simeq Z^{b_i}\oplus$(torsion). 
Let $P^{ab}$ be the abelianization of $P$. From above it follows that
$Z_Q(P)=Z_Q(P^{ab})$. Moreover, there is a 2-deformation
$$P^{ab}\to (\vee^{b_1}S^1)\vee T\vee (\vee^{b_2}S^2),$$
where $T$ has torsion homology and $\vee$ denotes 1-point union. 
Since $Z_{Q}$ is multiplicative under
1-point unions the corollary reduces to showing $Z_{Q}(S^1)=1$ and 
$Z_{Q}(S^2)=0$. The second
computation can be done formally. Notice these are dual in the sense 
that if $P=S^1$, which as a
presentation  has a single generator and no relations, then the dual has no 
generators and one
relation, so is $S^2$. The relation between dual presentations gives
$$X^{2}\,\hat Z(W_{S^1}) =\hat Z(W_{S^2}).$$
Hence $Z_{Q}(S^2)=0$.
More directly these are both defined in terms of the link invariant 
$Z(L)$ where $L$ is a single
unlinked circle. It is a simple computation that $Z(L)=X^2$. The 
invariant of $S^1$ is obtained by
dividing this by $X^2$ and reducing mod $p$, so $Z_{Q}(S^1)=1$ as required.
\end{ssec}

\begin{ssec}
The above corollary is not true in the case when the Euler 
characteristic of the complex is smaller or equal to 0. The simplest 
example comes from comparing the invariant of the 2-complexes 
$P=<x,y|\; xyx^{-1}y^{-1}>$ 
and $P'=<x,y|\; \emptyset >$, where $\emptyset$ denotes the trivial 
relation. The invariant of $P'$ is 0 (the link diagram of 
$W_{P'}$ consists 
of one undotted unknot and two dotted unknots, all disjoint from each 
other). The invariant of $P$, instead, is equal to the 
cardinality  $|S|$ of the set of simple objects $S$. 
The last asertion is proven in \cite{B2} but 
can also be seen using the present framework. In fact figure 
\ref{ab.gr} shows that
$$
\hat{Z}(W_P)=\sum_{b\in S}v^{t(b)}.
$$
Hence $Z_Q(P)=\phi_{p}(\hat{Z}(W_P))=|S|$.

\begin{figure}[h]
\setlength{\unitlength}{1cm}
\begin{center}
\begin{picture}(15,4.5)
\put(0,0){\fig{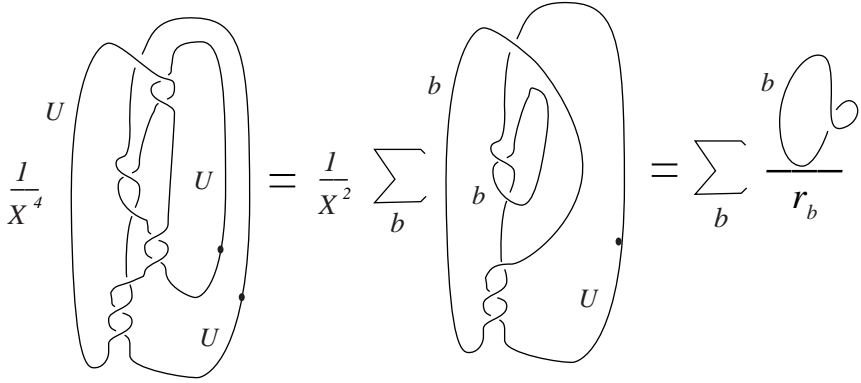}}
\end{picture}
\end{center}
\caption{Evaluation of $\hat{Z}(W_P)$ for $P=<x,y|\; xyx^{-1}y^{-1}>$.}
\label{ab.gr}
\end{figure}

\end{ssec}
\vspace{1cm}

\newsec{4}{Evaluation of  the $SL(2)$  2-complex invariant.}

\begin{ssec}
Here we give the proof of Proposition \ref{propa}, explicitly 
describing the 2-complex invariant for a
simple class of categories. According to \ref{corra} when the Euler 
characteristic is greater or equal to 1, the invariant 
depends only on the homology of the
complex. Given such 2-complex there is one with the same homology which 
is a 1-point union of copies of
$S^1$, $S^2$, and cyclic presentations $<x\mid x^n>$ for $n$ a power 
of a prime. The invariant is
multiplicative under 1-point unions and we have already found the 
values for $S^1$ and $S^2$. The
proposition is therefore reduced to proving the cyclic case:

\begin{lem} Let \A\ be the category obtained from
  class-0 $SL(2)$ representations at a $p^{th}$ root of unity, $p\geq5$. Then
$$Z_{Q}(<x\mid x^n>)=\left\{\begin{array}{cc} 0  &\mbox{if $p$ divides $n$};\\
                                n^{-2}\in Z/pZ , &\mbox{ otherwise}.
\end{array} \right.
$$
\end{lem}

   The category has simple objects  the  simple
representations of the quantum enveloping  algebra of $SL(2)$ with highest
weights $w$ satisfying
\begin{itemize}
\item[(i)] $w$ is even;
\item[(ii)] $0\leq w\leq p-3$.
\end{itemize}
The formulas here are taken from \cite{B1}, where the
definition of quantum enveloping algebra is the one of Lusztig, introduced
in \cite{L:book}. This is slightly different  from the definition
used in \cite{RT} and the literature
following it, and this leads to slight differences in the formulas.
\end{ssec}

\begin{ssec} Let $w=2z$ be a weight satisfying the conditions above. 
Then for the rank and the twist
of the representation with highest weight $w$, we have:
$$
  r_w=[2z+1],\quad\mbox{and }\quad
  \theta _w=v^{-2z(z+1)},
$$
where $[n]=\frac{v^n-v^{-n}}{v^1-v^{-1}}$. $[n]$ is often called a 
``quantum integer.''

For completeness we  list  formulas from \cite{nt.book} which
are used  below. Let
$$
g_{1}=\sum_{z=0}^{p-1}v^{z^2}=\sqcap
_{k=1}^{\frac{p-1}{2}}(v^{2k-1}-v^{-2k+1})=(v-v^{-1})^{\frac{p-1}{2}}
\sqcap _{k=1}^{\frac{p-1}{2}}[2 k+1]
\in
(v-1)^{\frac{p-1}{2}}k_{p},
$$
denote the Gauss sum, and
$$
\left(\frac{a}{p}\right) =\left\{\begin{array}{l}
0\quad \mbox{ if $a=0$ \quad $mod\quad p$},\\
1\quad \mbox{ if $a$ is a square $ mod\quad p$},\\
-1\quad \mbox{  if $a$ is not a square $ mod\quad p$},
\end{array}\right.
$$
be the Legendre symbol. Then in $R$, $p=(-1)^{\frac{p-1}{2}}g_1^2=
(v-v^{-1})^{p-1}
\sqcap _{k=1}^{\frac{p-1}{2}}[2 k+1]^2$.
For the present categories, using {\it Mathematica} and
the formula above, we obtain that
$$
X^2=\frac{-p}{(v-v^{-1})^2}=\frac{(-1)^{\frac{p+1}{2}}}{(v-v^{-1})^2}g_1^2
\in (v-1)^{p-3}k_p.
$$
Note this shows the mod $p$ reduction vanishes,
  as mentioned above.
\end{ssec}

\begin{ssec}
According to \ref{thm}, the invariant of the cyclic 2-complex
is the mod $p$ reduction of a normalized link invariant:
$$
Z_{Q}(<x\mid x^n>)=\phi_p(\hat Z(L_{n}))=\phi_p(\frac{1}{X^{2}}Z(L_{n})),
$$
where the link $L_{n}$ is shown in
figure \ref{cyc.diagr} (a). This can be transformed by a handle slide 
to the link in figure
\ref{cyc.diagr} (b). This shows the link invariant is a product with 
factors corresponding
to the two components, namely $F(n)F(-n))$, where
$$
F(n)=\sum_{z=1}^{\frac{p-3}{2}}(r_{2z})^2(\theta_{2z})^{n}=
\sum_{z=1}^{\frac{p-3}{2}}[2z+1]^2v^{-2zn(z+1)}.
$$
If
  $n$ is  divisible
by $p$ then $F(n)=F(-n)=X^2$, so the 2-complex invariant is the mod 
$p$ reduction of $X^2$ and
therefore vanishes.

\begin{figure}[h]
\setlength{\unitlength}{1cm}
\begin{center}
\begin{picture}(8,6)
\put(0,0){\fig{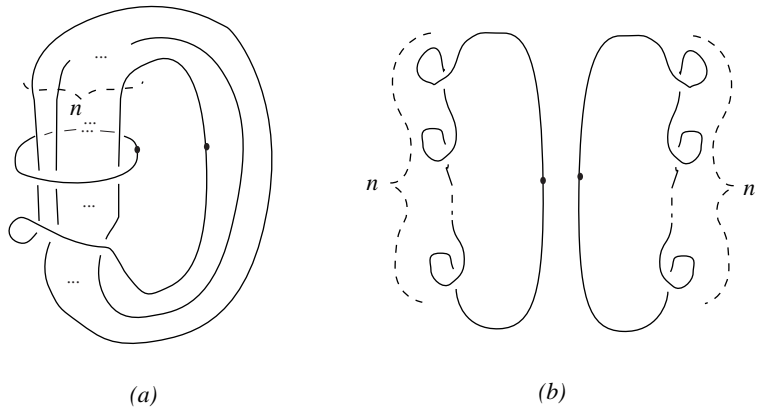}}
\end{picture}
\end{center}
\caption{The diagram corresponding to the cyclic group.}
\label{cyc.diagr}
\end{figure}

We now proceed  with  values of $n$
not divisible by $p$. First observe that the 3-manifold associated to 
the link in \ref{cyc.diagr} (b)
is the connected sum of two Lens spaces. Quantum invariants
of Lens spaces, as can be seen from the formula for $F(n)$,
  are reduced to calculating Gauss
sums, and have been studied for example in \cite{Stav}. Since the
category here is somewhat different the calculation needs to be
redone, but in a similar way we obtain
$$
F(n)=\left(\frac{-n/2}{p}\right)
\frac{g_{1}v^{\frac{n^2+1}{2n}}}{(v-v^{-1})}[\bar{n}],
$$
where
$\bar{n}$ denotes the inverse of $n$ in $Z/pZ$. Putting this in the 
expression for $Z_{Q}(<x\mid
x^n>)$ and reducing gives $\bar{n}^2$, as required.

\end{ssec}

\end{document}